\def \version {2014--06--28} 
\documentclass[12pt]{article}
\usepackage{graphpap}
\usepackage{stmaryrd}
\usepackage{amssymb}
\usepackage{amsmath}
\usepackage{algorithm}
\usepackage[noend]{algorithmic}

\def \qed {\hfill $\boxempty$}

\newtheorem{lem}{Lemma}
\newtheorem{Theorem}{Theorem}
\newtheorem{defi}{Definition}
\newtheorem{crl}{Corollary}
\newtheorem{prp}{Proposition}
\newtheorem{rmk}{Remark}
\newtheorem{xmp}{Example}
\newtheorem{clm}{Claim}
\newtheorem{op}{Problem}
\newtheorem{con}{Conjecture}

\def \bp {\begin{prp} \ }
\def \ep {\end{prp}}

\def \bc {\begin{crl} \ }
\def \ec {\end{crl}}

\def \bcon {\begin{con} \ }
\def \econ {\end{con}}

\def \thm {\begin{Theorem} \ }
\def \ethm {\end{Theorem}}

\def \bl {\begin{lem} \ }
\def \el {\end{lem}}

\def \bd {\begin{defi} \ \rm }
\def \ed {\end{defi}}

\def \brm {\begin{rmk} \ }
\def \erm {\end{rmk}}

\def \bxm {\begin{xmp} \ \rm }
\def \exm {\end{xmp}}

\def \bcm {\begin{clm} \ }
\def \ecm {\end{clm}}

\def \bop {\begin{op} \ }
\def \eop {\end{op}}

\def \nmr {\begin{enumerate}}
\def \enmr {\end{enumerate}}

\def \tmz {\begin{itemize}}
\def \etmz {\end{itemize}}

\def \nin {\noindent}
\def \bsk {\bigskip}

\def \pf {\nin{\bf Proof } \ }
\def \qed {\hfill $\Box$}

\begin{document}

\title {On the game domination number of graphs with given minimum degree\footnote{Research supported by the European Union and
Hungary and co-financed
 by the European Social Fund through the project T\'AMOP-4.2.2.C-11/1/KONV-2012-0004 - National Research Center
 for Development and Market Introduction of Advanced Information and Communication Technologies.}}

\author{Csilla Bujt\'{a}s\footnote{Email: bujtas@dcs.uni-pannon.hu }
\\
\normalsize
Department of Computer Science and Systems Technology \\
 \normalsize University of Pannonia\\
    \normalsize Veszpr\'{e}m,   Hungary
}
\date{\footnotesize Latest update \vspace{-6mm} on \version}
\maketitle

\begin{abstract}
 In the domination game, introduced by Bre\v{s}ar,  Klav\v{z}ar and   Rall in 2010,
 Dominator and Staller alternately select a vertex of a graph
 $G$. A move is legal if the selected vertex $v$  dominates at least one
 new vertex -- that is, if we have a $u\in N[v]$ for which no vertex
 from $N[u]$ was chosen up to this point of the game.
 The game ends when no more legal moves can be made, and its length
 equals the number of vertices selected.
 The goal of Dominator is to minimize whilst that of Staller is to maximize the
 length of the game. The game domination number
   $\gamma_g(G)$ of $G$ is the length of the domination game in which Dominator
   starts and both players play optimally.
  In this paper we establish an upper bound on $\gamma_g(G)$
  in terms of the minimum degree $\delta$ and the order $n$ of $G$.
  Our main result states that for every $\delta \ge 4$,
      $$\gamma_g(G)\le \frac{30\delta^4-56\delta^3-258\delta^2+708\delta-432}{90\delta^4-390\delta^3+348\delta^2+348\delta-432}\; n.$$
      Particularly,  $\gamma_g(G)
  < 0.5139\; n$ holds for every graph of minimum degree 4, and
    $\gamma_g(G) < 0.4803\; n$  if the minimum
    degree is greater than 4.
    Additionally, we prove that $\gamma_g(G) < 0.5574\; n$ if
    $\delta=3$.
\bigskip


\bigskip

\nin \textbf{AMS 2000 Subject Classification:}
 05C57,
 91A43,
 05C69

\end{abstract}

\vfil

\section{Introduction}

 In this note, our subject is the  domination game
   introduced by Bre\v{s}ar,  Klav\v{z}ar and   Rall in
 \cite{BKR-SIAM}.

 \subsection{Basic definitions}
   For a simple undirected graph $G=(V,E)$ and for a vertex $v\in V$,  the {\em  open
  neighborhood} of $v$  is  $N_G(v)=\{u : uv \in E\}$, while
  its {\em closed neighborhood} is   $N_G[v]=N_G(v) \cup \{v\}$.
  Then the {\em degree} $d_G(v)$  of $v$ is just
  $|N_G(v)|$ and the minimum degree  $\min \{d_G(v): v \in
  V\}$  is denoted by $\delta(G)$.
  As usual, we will write $N(v)$, $N[v]$ and $d(v)$ for $N_G(v)$, $N_G[v]$ and
  $d_G(v)$, respectively, if $G$ is clear from the context.

  Each vertex dominates itself and its neighbors, moreover
   a set $S\subseteq V$ dominates exactly those  vertices which are contained in
   $N[S]=\bigcup_{v\in S}N[v]$. A vertex set $D \subseteq V$ is
   called a {\em dominating set} of $G$ if  $N[D]=V$.
   The smallest cardinality of a dominating set is the {\em
   domination number}\/ $\gamma(G)$ of $G$.
  \bsk

  The {\it domination game}, introduced by Bre\v{s}ar,  Klav\v{z}ar and   Rall
 \cite{BKR-SIAM}, is played on a  simple undirected graph $G=(V,E)$ by  two players, named  Dominator and
 Staller, respectively. They  take turns choosing a vertex from $V$ such
 that a vertex $v$ can be chosen only if it dominates
 at least one new vertex -- that is, if we have a $u\in N[v]$ for which no vertex
 from $N[u]$ was selected up to this turn of the game.
 The game is over when no more legal
 moves can be made; equivalently, when the set $D$ of vertices chosen by the two players becomes a dominating set of
 $G$.
 The aim of Dominator is to finish the game as soon as possible, while that of Staller is to
  delay the end of the game.
 The {\it game domination number}\/ $\gamma_g(G)$  is the
 number of turns in the game
  when the first turn is Dominator's move  and both  players play
  optimally.
 Analogously, the {\it Staller-start game domination number}\/ $\gamma_g'(G)$ is the
 length of the game  when Staller begins   and the players   play
 optimally.

 \bsk
 \subsection{Results}

 Although the subject is quite new, lots of interesting results have been
 obtained on the domination game (see \cite{BDKK, BKKR, BKR-SIAM,
 BKR, JH, CS1, DKR, KWZ, Kos}). Note that also the total version of
 the domination game was introduced  \cite{HKR} and studied
 \cite{HKR2} recently.
 \bsk

  Concerning our present work, the bounds proved for the game domination number $\gamma_g(G)$
  are the most important preliminaries.
  The following fact was verified  in \cite{BKR-SIAM} and
   \cite{KWZ} as well.
  \begin{equation} \label{gamma}
  \gamma(G) \le \gamma_g(G) \le 2\gamma(G)-1
  \end{equation}
  Upper bounds in terms of the order were inspired by   the  following ``3/5-conjecture''
  raised by Kinnersley,  West and Zamani  \cite{KWZ}.
  \begin{con}
  \label{conj1}
  If $G$ is an isolate-free graph of order $n$, then
  $\gamma_g(G) \le 3n/5$
  holds.
  \end{con}
  Conjecture \ref{conj1} has been proved for the following graph
  classes:
  \tmz
  \item   for trees of order $n \le 20$ (Bre\v{s}ar, Klav\v{z}ar, Ko\v{s}mrlj and
  Rall  \cite{BKKR});
  \item for caterpillars -- that is, for trees in which the non-leaf
  vertices induce a path (Kinnersley,  West and Zamani \cite{KWZ});
  \item for trees in which no two leaves are at distance four apart
  (Bujt\'as \cite{JH, CS1}).
  \etmz
  Moreover, 
  in a manuscript in preparation,  Henning and Kinnersley prove
  Conjecture~\ref{conj1} for graphs of minimum degree at least 2 \cite{HK}.

  \bsk

  On the other hand, upper bounds weaker than $3n/5$ were obtained for
  some wider graph classes. For trees, the inequality
    $\gamma_g(G) \le  7n/11$  was established by
   Kinnersley,  West and Zamani in \cite{KWZ} and  it was recently
   improved to $\gamma_g(G) \le 5n/8$ by the present author in \cite{CS1}.
   For the most general case,
   Kinnersley,  West and Zamani proved \cite{KWZ} that the game
   domination number of any isolate-free graph $G$ of order $n$ satisfies
   $\gamma_g(G)\le \lceil 7n/10 \rceil$.
   In Section \ref{sect2} we improve this upper bound by
   establishing the following claim.
   \bp \label{2/3}
    For any isolate-free graph $G$ of order $n$,
    $$\gamma_g(G)\le \frac{2n}{3} \qquad \mbox{and} \qquad \gamma'_g(G)\le
    \frac{2n}{3}.$$
   \ep
   In fact, in a manuscript under preparation \cite{CS2} we will
   prove the stronger inequality $\gamma_g(G)\le 0.64 n$, but the
   proof of Proposition \ref{2/3} may be of interest because of its
   simplicity and gives  illustration for the proof technique applied
   in the later sections.

   \bsk
   One of our main results gives an upper bound smaller than $0.5574n$ on
  the game domination number of  graphs with minimum degree   3.

   \thm \label{deg3}
   For any graph $G$ of order  $n$ and with minimum degree 3,
   $$\gamma_g(G)\le  \frac{34 n}{61}  \qquad \mbox{and} \qquad
   \gamma'_g(G)\le  \frac{34 n-27}{61}.$$
   \ethm

 For graphs all of whose vertices are of degree greater than 3, we
 prove an upper bound   in terms of the order
   and the minimum degree.
   \thm \label{main}
    If  $G$ is a graph on $n$ vertices and its minimum degree is
    $\delta(G)\ge d \ge 4$, then
        \begin{eqnarray*}
        \gamma_g(G) &\le & \frac{30d^4-56d^3-258d^2+708d-432}{90d^4-390d^3+348d^2+348d-432}\;
        n.
        \end{eqnarray*}
   \ethm

   As the coefficient in this upper bound equals $37/72 <0.5139$ for $d=4$, and equals $2102/4377<0.4803$ for $d=5$,
    the following immediate consequences are obtained.

   \bc \label{deg4}
   \tmz
   \item[$(i)$] For any   graph $G$ of order $n$ and with minimum degree
   $\delta(G)= 4$, the inequality
    $\gamma_g(G)\le 37n/72$
    holds.
    \item[$(ii)$] For any   graph $G$ of order $n$ and with minimum degree
   $\delta(G)\ge 5$, the inequality
    $\gamma_g(G)\le 2102n/4377$
    holds.
    \etmz
   \ec
   Particularly, these statements  show that the coefficient $3/5$  in Conjecture \ref{conj1}
   can be significantly improved if only those graphs with
    $\delta(G) \ge 4$ are considered.

   On the other hand, note that Theorem~\ref{deg3} and Theorem~\ref{main} establish new
   results
   only for $3\le \delta(G) \le 21$. Although it was not mentioned
   in  the earlier papers,  the upper bound in (\ref{gamma})
   together with the well-known theorem (see e.g., \cite{AS})
   \begin{equation*}
   \gamma(G)\le \frac{1+ \ln (\delta +1)}{\delta+1}\; n
   \end{equation*}
    clearly  yields
   \begin{equation}\label{A-S}
   \gamma_g(G) < 2\cdot \frac{1+ \ln (\delta +1)}{\delta+1}\; n
   \end{equation}
   for each $\delta \ge 2$.
   For integers $3\le \delta(G)=d \le 21$, it is easy to check that our
   bound is better than the above one in (\ref{A-S}).

   \bsk

  Our proof technique is based on a value assignment to the
  vertices where the value of a vertex depends on
   its current status  in the game.  We will consider a greedy strategy of Dominator, where
  the greediness is meant concerning the decrease in the values.
   Our main goal is to estimate the
   average decrease in a turn  achieved under this assumption.
   We have been introduced this type of approach in the conference paper
   \cite{JH} and in the paper \cite{CS1}. The frame of this technique and the basic
   observations are contained here in Section~\ref{sect2}.
   Then, in Section~\ref{sect3} and Section~\ref{sect4} we specify
   the details and prove our Theorem~\ref{deg3} and
   Theorem~\ref{main} respectively. In the last section we make some additional
   notes concerning the Staller-start version of the game.

\bsk



\section{Preliminaries} \label{sect2}

 Here we introduce the notion of the residual graph, define the color assignment to the
 vertices and
 give a general determination for the phases of the game. Then, we
 take some simple observations which will be used in the
 later sections.

 \paragraph{Colors}
 Consider any moment of the process of a domination game on
 the graph $G^*=(V,E)$, and denote by $D$ the set of  vertices
 chosen up to this point of the game. As it was introduced in \cite{JH} and \cite{CS1}, we distinguish between the
 following three types of vertices.
 \tmz
 \item A vertex $v\in V$ is \emph{white} if $v \notin N[D]$.
 \item A vertex $v\in V$ is \emph{blue} if $v \in N[D]$ but $N[v] \nsubseteq N[D]$.
 \item A vertex $v\in V$ is \emph{red} if  $N[v] \subseteq N[D]$.
 \etmz

 \paragraph{Residual graph} Clearly, a  red vertex $v$ and all its
 neighbors are already dominated in the game. Hence the choice of $v$
 would not be  a legal move in the later turns and further, the
 status of $v$ remains red. So, red vertices do not influence the
 continuation of the game and they can be deleted. Similarly, edges
 connecting two blue vertices can be omitted too. This graph,
 obtained after the deletion of red vertices and edges between two
 blue vertices, is  called \emph{residual graph}, as it was
 introduced in \cite{KWZ}. At any point of the game, the set of
 vertices chosen up to this point is denoted by $D$ and the residual
 graph is denoted by $G$. When it is needed, we use the more precise notations $D_i$ and $G_i$
 for the current $D$ and $G$ just before
 the $i$th turn. 

 \paragraph{Phases of the game}
 The phases will be defined for the Dominator-start game that is,
 for each odd integer $j$ the $j$th turn belongs to Dominator. The
 Staller-start version will be treated later by introducing a Phase
 0 for the starting turn.

 In our proofs, nonnegative values $p(v)$ are assigned to the
 vertices, and the value $p(G)$ of the residual graph is just the
 sum of the values of the vertices.
 Also, we assume that Dominator always chooses greedily. More precisely, for  each odd $j$, in the $j$th turn  he plays a
 vertex which results the possible maximum  $p(G_j)-p(G_{j+1})$.
 This difference is called the decrease in the value of $p(G)$ and
 also referred to as the gain of the player.

 \bd \label{def-phases}
 Let $(C1),\dots (C\ell)$ be conditions all of which relating to the
 $j$th turn of the game where $j$ is odd.
 Then, for each $i=1,\dots \ell$, Phase $i$ of the game is defined as
 follows.
  \tmz
  \item[$(i)$] Phase 1 begins with the first turn of the game.
 \item[$(ii)$] If Phase $i$ begins with the $b_i$th turn, it is
 continued as long as $(Ci)$ is satisfied in each turn of Dominator.
 That is, Phase $i$ ends right after the $e_i$th turn where $e_i$ is
 the smallest even integer with $b_i <e_i$ for which $(Ci)$ is not
 satisfied in the $(e_i+1)$st turn.
  \item[$(iii)$] If Phase $i$ ends after the $e_i$th turn but the
  game is not over yet, then the $(e_i+1)$st turn is the beginning
  of Phase $i'$, where $i'$ is the smallest integer with $i<i'$ such
  that $(Ci')$ is fulfilled in the $(e_i+1)$st turn.
  \item[$(iv)$] If Phase $i$ is followed by Phase $i'$ and $i+2 \le
  i'$ holds, we say that Phases $i+1, \dots i'-1$ are skipped;
  moreover, their starting and end points are interpreted to be
  the same as the end of Phase $i$.
 \etmz
 \ed

\paragraph{Further notations}
The colors white, blue and red will be often abbreviated to W, B and
R, respectively. For example, a B-neighbor is a blue neighbor, and
the notation $v$: W$\rightarrow$B/R means that vertex $v$ changed
from white to either blue or red in the turn considered. Moreover,
  $d^W(v)$ and $d^B(v)$ stand for the number of W-neighbors
 and  B-neighbors of $v$, respectively.

 \bigskip

 We  cite
 the following observations (in a slightly modified form) from \cite{CS1}:
 \bl \label{lemma-basic}
 The following statements are true for every residual graph $G$ in a domination game started on
 $G^*$.
 \tmz
 \item[$(i)$]
  If $v$ is a white vertex in  $G$, then $v$ has
  the same neighborhood in $G$ as it had in $G^*$. Thus, $d_G(v)=d_{G^*}(v)$ holds for every W-vertex of
    $G$ and moreover, $d^W_G(v)+d^B_G(v)=d_{G^*}(v)$.
\item[$(ii)$]
  If $v$ is a blue vertex in $G$, then $v$ has only white neighbors
  and definitely has at least one. That is, $d^W_G(v)=d_G(v)\ge 1$ and $d^B_G(v)=0$ if $v$
  is a B-vertex in $G$.
 \etmz
 \el

 At the end of this section, we provide a simple example for applying the tools introduced above.
 We prove Proposition~\ref{2/3}, which states that for any isolate-free graph $G$ of order $n$,
 $$\gamma_g(G)\le \frac{2n}{3} \qquad \mbox{and} \qquad \gamma'_g(G)\le
    \frac{2n}{3}$$
    hold.

 \bsk

 \noindent \textbf{Proof of Proposition~\ref{2/3}.}
 First, we consider the Dominator-start game on $G^*=(V,E)$, which is a simple graph without isolated vertices.
 In every residual graph $G$, let the value $p(v)$ of a vertex $v$ be equal to 2, 1 and 0,
 when $v$ is white, blue and red, respectively. Hence, we start with
 $p(G^*)=2n$ and assume that Dominator always selects a vertex which results in a maximum decrease in $p(G)$.
  The game is divided into two phases, which are
 determined due to Definition~\ref{def-phases} with the following
 conditions:
 \tmz
 \item[$(C1)$] Dominator gets at least 4 points.
 \item[$(C2)$] Dominator gets at least 1 point.
 \etmz
 \paragraph{Phase 1.} If Staller  selects  a
 W-vertex, then it becomes red and causes at least 2-point decrease
 in the value of the residual graph. In the other case, Staller
 selects a B-vertex $v$ which has a W-neighbor $u$. Then, the
 changes $v$: B$\rightarrow$R and $u$: W$\rightarrow$B/R together result
 in a decrease of at least $1+1=2$. Hence, in each of his turns
 Staller gets   at least 2 points. By condition $(C1)$, Dominator
 always gets at least 4 points. As Dominator begins the phase, the
 average decrease in $p(G)$  must be at least 3 in a turn.
 \paragraph{Phase 2.} When Phase 2 starts, Dominator cannot seize 4
 or more points by playing any vertex of $G_j$. This implies the
 following properties of the residual graph:
 \tmz
 \item For every W-vertex $v$, $d^W(v) \le 1$.\\
  Indeed, if $v$ had two W-neighbors $u_1$ and $u_2$, then Dominator
 could choose $v$ and the changes
 $v$: W$\rightarrow$R and $u_1, u_2$: W$\rightarrow$B/R would give a
 gain of at least $2+2\cdot 1=4$ points, which is a contradiction.
 \item For every W-vertex $v$, $d^W(v) = 0$. \\
  We have seen that $d^W(v) \le 1$. Now, assuming two W-neighbors $v$
 and $u$, the choice of $v$ would result in the changes
 $v, u$: W$\rightarrow$R, which give a gain of at least 4 points to
  the player. This is a contradiction again.
 \item For every B-vertex $v$, $d(v) = 1$. \\
 By Lemma~\ref{lemma-basic}$(ii)$, $d(v)\ge 1$. Now, assume that $v$
 has  two different W-neighbors $u_1$ and $u_2$. As we have shown,
 $d^W(u_1) =d^W(u_2) = 0$ must hold and consequently, if Dominator plays
 $v$, then both $u_1$ and $u_2$ turn to  red. This gives a gain
 of at least $1+2\cdot 2=5$ points, which cannot be the case at the
 endpoint of Phase 1.
 \item Each component of $G_j$  is a $P_2$ and contains exactly
 one white and  one blue vertex.\\
 By the claims above, each component is a star with a white center
 and blue leaves. If it contained at least two leaves then Dominator
 could play the center and get at least 4 points.
 \etmz

Therefore, at the beginning of Phase 2  the residual graph consists
of components of order 2. As follows, in each turn an entire
component becomes red and $p(G)$ decreases by exactly 3 points.

In the game, the value of the residual graph decreased from $2n$ to
zero, and the average decrease in a turn was proved to be at least
3. Consequently, the number of turns required is not greater than
$2n/3$, which proves $\gamma_g(G^*) \le 2n/3$. \bsk

 If Staller starts the game, his first move definitely decreases
 $p(G)$ by at least 3 points as there are no isolated vertices. Then, the
 game is continued as in the Dominator-start game, and the average
 decrease remains at least 3 points. Thus, $\gamma_g'(G^*) \le 2n/3$ holds.

\qed \bsk

 \section{Graphs of minimum degree 3} \label{sect3}

 In this section we prove the upper bound stated on the game domination
 number of graphs with minimum degree 3. Also, this proof serves as
 an introduction to the details of the idea applied in the next section to prove our main
 theorem.
 \bsk

 \noindent \textbf{Proof of Theorem~\ref{deg3}.}
  We consider a graph $G=(V,E)$ of minimum degree 3 and define the
   value assignments of types A1.1, A1.2 and A1.3 as they are given in
   Table~\ref{table-deg3}.
\begin{table}[!h]
 \nin \caption{Value assignments used in the proof of Theorem~\ref{deg3}}  \label{table-deg3}
\begin{center} \begin{tabular}{|c|c|c|c|c|}\hline
   Abbreviation
  &Type of the vertex
 &Value in A1.1
 &Value in A1.2
 &Value in A1.3
  \\ \hline \hline
  W
 &white vertex
 &$34$
 &$34$
 &$34$
  \\ \hline 
  B$_3$
 &blue vertex of degree at least 3
 &16
 &16
 &---
  \\ \hline
  B$_2$
 &blue vertex of degree 2
 &16
 &13
 &13
  \\ \hline
  B$_1$
 &blue vertex of degree 1
 &16
 &10
 &9
  \\ \hline
R
 &red vertex
 &0
 &0
 &0
  \\ \hline

\end{tabular}
 \end{center}
  \end{table}

\nin Hence, the game starts with $p(G^*)=34 n$. First, assume that
Dominator begins the game and determine Phases 1-4 due to
Definition~\ref{def-phases} with the following specified conditions:
  \tmz
  \item[$(C1)$] Dominator gets at least 88 points due to the
  assignment A1.1.
  \item[$(C2)$] Dominator gets at least 91 points due to the
  assignment A1.2.
  \item[$(C3)$] Dominator gets at least 84 points due to the
  assignment A1.3.
  \item[$(C4)$] Dominator gets at least 1 point  due to the
  assignment A1.3.
  \etmz

  \paragraph{Phase 1.} Here, we apply the value assignment A1.1.
  In each of his turns, Staller either selects
  a white vertex and gets at least 34 points; or he plays a blue
  vertex $v$ which has a white neighbor $u$ and then the color changes
  $v$: B$\rightarrow$R and $u$: W$\rightarrow$B/R give  at least
  16+18=34 points. By condition $(C1)$, each move of Dominator
  yields a gain of at least 88 points. As Dominator begins the game,
  we have the following estimation on the average decrease of $p(G)$
  in a turn.
  \bl \label{ph1}
  In  Phase 1, the average decrease of $p(G)$ in a turn  is at
  least $61$ points.
  \el
  At the end of Phase 1 we have some structural properties which
  remain valid in the continuation of the game.
  \bl \label{str-ph1}
   After the end of Phase 1, throughout the game, each white vertex
   has at most 2 white neighbors, and each blue vertex has at most 3
   white neighbors.
   \el
  By definition, at the end of the first phase Dominator has no possibility to
  seize 88 or more points by playing a vertex of the residual graph
  $G$.
   Now, assume that there exists a W-vertex $v$ with three W-neighbors
  $u_1$, $u_2$ and $u_3$ in $G$. Then, Dominator could choose
  $v$ and the color changes
   $v$: W$\rightarrow$R and $u_1, u_2, u_3$: W$\rightarrow$B/R
   would   decrease $p(G)$ by  at least $34+3 \cdot 18= 88$ points,
   which is a contradiction.
   Similarly, if there exists a B-vertex $v$ with at least four
   W-neighbors, then Dominator could get at least $16+4\cdot 18=88$
   points by playing $v$, which is a contradiction again.

   In the continuation of the game, new white vertices cannot arise,
   moreover a new blue vertex may arise only by the color change
   W$\rightarrow$B. This implies that the stated properties remain
   valid throughout all the later phases. \qed

  \paragraph{Phase 2.} At the beginning of this phase we change to
  assignment A1.2. Clearly,  the
  values  of the vertices do not increase. As no  blue vertex  has
  a degree greater then 3, we observe that each B-vertex $v$ has value
  $p(v)=7+ 3d(v)$.
  Further, assignment A1.2
  ensures that when a vertex $v$ is played, the value  of every blue
  vertex from $N[N[v]]$ is decreased.
  \bl \label{gain-2}
  The following statements are true in Phase 2.
  \tmz
   \item[$(i)$] If a W-vertex\/ $v$ with W-degree\/ $d^W(v)$ is played,
    then\/ $p(G)$ decreases by at least\/ $43+24 d^W(v)$ points.
    \item[$(ii)$] If a B-vertex\/ $v$ with degree\/ $d(v)$ is played,
    then\/ $p(G)$ decreases by at least\/ $7+24 d(u)$ points.
    \item[$(iii)$] In each turn\/ $p(G)$ decreases by at least\/ $31$ points
  \etmz
  \el
  \pf When the degree of a B-vertex is decreased by $x$, its value
  decreases by at least $3x$, no matter whether the change is of
  type  B$_i \rightarrow$B$_{i-x}$ or B$_x \rightarrow$R. Thus, if a vertex $v$ is played,
  the
  sum of the values of B-vertices contained in $N[N[v]]\setminus \{v\}$ is
  decreased by at least
  $$3\sum_{u \in N[v]} d^B(u)\ge  3\sum_{u \in N[v]} (3-d^W(u))$$
  if $v$ is
  white, and by at least
  $$3\sum_{u \in N[v]} (d^B(u)-1) \ge 3\sum_{u \in N[v]} (2-d^W(u))$$
   if
  $v$ is blue.
  \bsk

  First, assume that the played vertex $v$ is white,   $d^W(v)=k$
  and the W-neighbors of $v$
   are $u_1, \dots  u_k$. For each $1\le i \le k$, the W-vertex
   $u_i$ becomes either a B-vertex of degree at most $d^W(u_i)-1$ or an
   R-vertex. As $1 \le d^W(u_i) \le 2$, $p(u_i)$ decreases by at least
   $34-(7+3d^W(u_i)-3)=30-3d^W(u_i)$
   in either case.
   Then, the decrease in $p(G)$ is not smaller than
   $$34+ \sum_{i=1}^k (30-3d^W(u_i)) + 3(3-k)+3\sum_{i=1}^k
   (3-d^W(u_i))=43+36k-6\sum_{i=1}^k d^W(u_i) \ge 43+24 k,
   $$
   where $0 \le k \le 3$ must hold. This establishes statement $(i)$.
   \bsk

   In the other case, $v$ is blue with $d(v)=k$ and its W-neighbors
   are $u_1, \dots  u_k$. As $v$ has only white neighbors and definitely has  at least
   one and no more than 3,
   $1 \le k \le 3$ holds; moreover, $0 \le d^W(u_i) \le 2$ is true
   for all $1\le i \le k$.
   When $v$ is played, $u_i$ becomes red if $d^W(u_i)=0$, otherwise
   it will be  a blue vertex of degree
   at most $d^W(u_i)$. Therefore, the decrease in $p(u_i)$ is at
   least $34-(7+3d^W(u_i))=27-3d^W(u_i)$ and that in $p(v)$ is
   exactly $7+3k$.
   Then, the sum of the decreases cannot be smaller than
     $$7+3k+ \sum_{i=1}^k (27-3d^W(u_i)) +3\sum_{i=1}^k
   (2-d^W(u_i))=7+36k-6\sum_{i=1}^k d^W(u_i) \ge 7+24 k
   $$
   as stated in $(ii)$.

   \bsk

   To prove $(iii)$, it suffices to  consider the minimum of
      $43+24 k$ in case $(i)$, which is 43; and that of $7+24 k$ in
      case $(ii)$, which is 31 because of the condition $k \ge 1$.
      \qed

      \bsk

      By Lemma~\ref{gain-2}$(iii)$,
      Staller gets at least 31 points, and by Condition $(C2)$, Dominator
        gets at least 91 points in each of their turns. Hence,   we have the
        following estimation.
      \bl \label{ph2}
  In   Phase 2, the average decrease of $p(G)$ in a turn is at
  least $61$ points.
  \el
   As shown by the next lemma,  the W-degrees are more strictly
   bounded from the end of Phase 2 than earlier.
   \bl \label{str-ph2}
   After the end of Phase 2, throughout the game, each white vertex
   has at most 1 white neighbors, and each blue vertex has at most 2
   white neighbors.
   \el
   \pf By condition $(C2)$, at the end of Phase 2 Dominator can
   seize only less than 91 points by choosing any vertex of $G$.
   By Lemma \ref{gain-2}$(i)$, the selection of a W-vertex $v$ with
   $d^W(v)=2$ causes a decrease of at least $43+24\cdot 2=91$ points
   in $p(G)$. Hence, each W-vertex has either zero or exactly one
   W-neighbor.

   Now, assume that $v$ is  a B-vertex   with three W-neighbors, say $u_1$, $u_2$ and
   $u_3$. We have already seen that the inequalities  $0 \le d^W(u_i) \le 1$ hold for
   $i=1,2,3$. Then, as it was shown in the proof of Lemma~\ref{gain-2}$(ii)$, the choice of $v$ would decrease   $p(G)$ by at
   least
   $$7+36\cdot 3-6\sum_{i=1}^3 d^W(u_i)\ge 97,$$
   which is a contradiction.
   \qed

   \paragraph{Phase 3.} The phase starts with changing the value
   assignment A1.2 to A1.3. By Lemma~\ref{str-ph2}, there are no B-vertices of degree 3 or
   higher, moreover we observe that the change to A1.3 cannot cause increase in
   the value of $G$. Also, one can easily check that the value of a
   B-vertex decreases by at least $4x$ points, if it loses $x$
   W-neighbors in a turn.

   \bl \label{gain-3}
  The following statements are true in Phase 3.
  \tmz
   \item[$(i)$] If a W-vertex\/ $v$  is played,
    then\/ $p(G)$ decreases by at least\/ $84$ points if $d^W(v)=1$,
    and\/
    $p(G)$ decreases by at least\/ $46$ points if\/ $d^W(v)=0$.
    \item[$(ii)$] If a B-vertex\/ $v$  is played,
    then\/ $p(G)$ decreases by at least\/ $67$ points if\/ $d(v)=2$,
    and\/
    $p(G)$ decreases by at least\/ $38$ points if\/ $d(v)=1$.
    \item[$(iii)$] In each turn\/ $p(G)$ decreases by at least\/ $38$ points
  \etmz
  \el
   \pf $(i)$ Consider a W-vertex $v$ whose only W-neighbor is $u$. By
    Lemma~\ref{str-ph2}, all the further neighbors of $v$ and $u$ are blue
   and there are at least four such neighbors.
    Hence, when $v$ is played, the color changes $v,u$: W$\rightarrow$R
   decrease $p(G)$ by 68 points, while the
  sum of the values of B-vertices contained in $N[\{v,u\}]$
  decreases by at least $4 \cdot 4=16$. Hence, the gain of the
  player is at least 84 points.
  In the other case, when  $v$ has no W-neighbors, it has at least three B-neighbors.
  Then, the change
  $v$: W$\rightarrow$R gives at least 34 points and additionally, the decrease in the degrees of the
   B-neighbors means  at least 12 points. This proves that
  $p(G)$ decreases by at least 46 points.

  $(ii)$ If the played vertex $v$ is blue and has exactly one white
  neighbor $u$, then the changes $v$: B$_1\rightarrow$R and
  $u$: W$\rightarrow$B$_1$/R cause a decrease of at least $9+25=34$ points
  in $p(G)$. Additionally, $u$ has at least one B-neighbor
  different from $v$, whose value is decreased by at least 4 points.
  Consequently, the total decrease is at least 38 points.
  Similarly, if $v$ is blue and has two W-neighbors $u_1$ and $u_2$, then the
  total decrease in $p(G)$ is at least $13+2\cdot 25+2\cdot 4=71$
  points.

  $(iii)$ As the four cases above cover all possible moves which can
  be made in Phase~3, $p(G)$ is decreased by at least 38 points in each
  turn. \qed
  \bsk

  As   consequences of Condition $(C3)$ and
  Lemma~\ref{gain-3}$(iii)$,
    Dominator gets at least 84 points and Staller gets at
  least 38 points in each of his turns. Hence, we have the
  desired average.

   \bl \label{ph3}
  In   Phase 3, the average decrease of $p(G)$ in a turn is at
  least $61$ points.
  \el
   When Dominator cannot get at least 84 points in a turn, the
   structure of the residual graph must be very simple.

  \bl \label{str-ph3}
   At the end of Phase 3,  each component of the residual graph is a star of order $k \ge 4$ with a
   white center and $k-1$ blue leaves.
      \el
   \pf Let $G_i$ be the residual graph obtained at the end of Phase 3.
   Due to Lemma~\ref{gain-3}$(i)$,  the presence of a W-vertex $v$ with $d^W(v)=1$
   provides an opportunity for Dominator to get at least 84 points.
   Then, the $i$th turn would belong to Phase 3, which is a
   contradiction. Consequently, in $G_i$ each W-vertex has only
   B-neighbors.

   Next, assume that we have a B-vertex $v$ which has two
   W-neighbors $u_1$ and $u_2$ in $G_i$. As we have seen,
   in $G_i$
   $d^W(u_1)=d^W(u_2)=0$ must hold  and moreover, both $u_1$ and
   $u_2$ have at least two B-neighbors.
   Therefore, if $v$ is selected by Dominator,
   the changes $v$: B$_2\rightarrow$R and
   $u_1, u_2$: W$\rightarrow$R  with the change in the
   values of B-neighbors, all together yield  at least
   $13+ 2\cdot 34+ 4 \cdot 4=97$ point decrease in $p(G_i)$, which
   is a contradiction again. Hence, each B-vertex has at most one
   W-neighbor.

   Since each W-vertex $v$ has the same degree in the residual
   graph $G_i$ as it had in $G^*$, it has at least three B-neighbors
   in $G_i$. In addition, each B-vertex is a leaf in $G_i$. This
   implies that at the
   end of Phase 3  every component is a star with the structure stated.
   \qed

   \paragraph{Phase 4.} By Lemma~\ref{str-ph3}, Phase 4 begins with
   star-components containing a white center and at least three blue
   leaves. Then, in each turn a component becomes completely red, no matter whether a white or a blue vertex is
   played. Thus, each move   decreases  the value of $G$ by at least $34+3
   \cdot 9= 61$ points.

   \bl \label{ph4}
  In   Phase 4, the average decrease of $p(G)$ in a turn is at
  least $61$ points.
  \el

  \bsk
  By Lemmas~\ref{ph1}, \ref{ph2}, \ref{ph3} and \ref{ph4}, if
  Dominator starts the game and he plays the prescribed greedy strategy, then for the number $t^*$ of turns
    $$\gamma_g(G^*) \le t^* \le \frac{34}{61}n$$
    holds.

  \bsk

  Finally, for the Staller-start version of the game we define
  Phase 0, which contains only the first turn and the values are counted due to A1.1. Observe that
  Staller's any choice results in
  at least $34+ 3\cdot 18=88$ point decrease in $p(G^*)$. Then, Phase 1
  might be skipped if $(C1)$ is not true for $G_1$,  but otherwise the game
  continues as in the Dominator-start version and   our lemmas
  remain valid. Therefore, by the 27 point overplus arising in Phase
  0, for $\gamma_g'(G^*)$ we obtain a slightly better bound,
  $$\gamma_g'(G^*) \le  \frac{34n-27}{61}.$$
  This completes the proof of Theorem~\ref{deg3}. \qed

\bsk

\section{Graphs with minimum degree greater than 3} \label{sect4}

Here we prove  Theorem  \ref{main} and Corollary \ref{deg4}.

\bsk

 \noindent \textbf{Proof of Theorem~\ref{main}.}
 First, we consider the Dominator-start game on a graph $G^*=(V,E)$ of order
$n$, whose minimum degree is $\delta(G^*) \ge d \ge 4$.

The proof and the game starts with the value assignment A2.1 to the
vertices as shown in Table \ref{table-deg4}.
 Later,  we use a more subtle distinction between the types of
blue vertices due to assignments A2.2, A2.3 and A2.4 (see
Table~\ref{table-deg4}.
  We will see that the value $p(G)$ of the residual graph cannot
increase when we change to an assignment with a higher index.
 \bsk
  We consider a graph $G^*=(V,E)$ with a minimum degree $\delta(G)\ge d\ge 4$ and define the
   value assignments of types A2.1, A2.2, A2.3 and A2.4 as they are given in
   Table~\ref{table-deg4}.
\begin{table}[!htb]
 \nin \caption{Value assignments used in the proof of Theorem~\ref{main}}  \label{table-deg4}
\begin{center} \begin{tabular}{|c|c|c|c|c|c|}\hline
   Abbreviation
  &Type of the vertex
 &A2.1
 &A2.2
 &A2.3
 &A2.4
  \\ \hline \hline
  W
 &white vertex
 &$a$
 &$a$
 &$a$
 &$a$
  \\ \hline 
  B$_4$
 &blue vertex of degree at least 4
 &$b$
 &$b$
 &---
 &---
  \\ \hline
  B$_3$
 &blue vertex of degree   3
 &$b$
 &$b-x_1$
 &$b-x_1$
 &---
  \\ \hline
  B$_2$
 &blue vertex of degree 2
 &$b$
 &$b-2x_1$
 &$b-x_1-x_2$
 &$b-x_1-x_2$
  \\ \hline
  B$_1$
 &blue vertex of degree 1
 &$b$
 &$b-3x_1$
 &$b-x_1-2x_2$
 &$b-x_1-x_2-x_3$
  \\ \hline
R
 &red vertex
 &0
 &0
 &0
 &0
  \\ \hline
\end{tabular}
 \end{center}
  \end{table}

  \nin The values of $a$, $b$, $x_1$, $x_2$, $x_3$ and $s$ are
  defined in terms of the parameter $d$. We aim to prove that $s$ is
  a lower bound on the average decrease of $p(G)$ in a turn, if
  Dominator follows the prescribed greedy strategy.
      \begin{eqnarray*}
  a&=& 30d^4-56d^3-258d^2+708d-432\\
  b&=& 111d^3-561d^2+888d-432\\
  x_1&=& 6d^3-19d^2+15d\\
  x_2&=& 15d^3-64d^2+65d\\
  x_3&=& 30d^3-144d^2+202d-72\\
  s&=& 90d^4-390d^3+348d^2+348d-432
  \end{eqnarray*}
  Concerning the values above and the change between assignments, we take the following observations.
  \bl \label{lemma-count}
  For every fixed integer $d\ge 4$:
  \tmz
   \item[$(i)$] $0 < x_1 < x_2 < x_3 < b-x_1-x_2-x_3<
    b<a$ \, and \, $x_3 <a-b$.
   \item[$(ii)$] For every $1\le i <j \le 4$ and every residual graph $G$,   $p(G)$ does not
   increase if the value assignment A$2$.$i$ is changed to A$2$.$j$
   (assuming that A$2$.$j$ is defined for $G$).
  \etmz
  \el
  \pf The proof of $(i)$ is based on a simple counting and estimation.
  Table~\ref{counting} shows the differences and their exact values  for
  $d=4,5,6$.
  The comparison of coefficients verifies our statements for $d\ge
  7$.

  \begin{table}[!htb]
 \nin \caption{Values of the differences for the proof of Lemma~\ref{lemma-count}}  \label{counting}
\begin{center} \begin{tabular}{|l|l|r|r|r|}\hline
     &
     &$d=4$
 &$d=5$
 &$d=6$
   \\ \hline \hline
  $x_1$
  &$6d^3-19d^2+15d$
 &140 &8408    &21312
  \\ \hline
  $x_2-x_1$
  &$9d^3   -45d^2+ 50d-72$
  &56  &250 &624
  \\ \hline
  $x_3-x_2$
  &$15d^3  -80d^2+ 137d-288$
  &156 &488 &1110
  \\ \hline
  $b-x_1-x_2-2x_3$
  &$30d^3  -190d^2+    404d$
  &208 &732 &1776
  \\ \hline
  $a-b$
  &$30d^4  -167d^3+    303d^2 -180d$
  &1120    &4550    &12636
  \\ \hline
  $a-b-x_3$
  &$30d^4  -197d^3+    447d^2 -382d+    72$
  &768 &3462    &10200
  \\ \hline

\end{tabular}
 \end{center}
  \end{table}

Once $(i)$ is proved, Table~\ref{table-deg4} shows that no vertex
  has greater value by A2.$j$ than by A2.$i$, whenever $j>i$ holds.
    \qed

    \bsk
  Note that later we will use further relations between
  $a$, $b$, $x_1$, $x_2$, $x_3$ and $s$    but these are equations,
  which can be verified by simple counting, so the details will be
  omitted. \bsk

The game is divided into five phases due to
Definition~\ref{def-phases} with the following five conditions:
  \tmz
  \item[$(C1)$] Dominator gets at least $5a-4b$ points due to the
  assignment A2.1.
  \item[$(C2)$] Dominator gets at least $4a-3b+(4d-6)x_1$ points due to the
  assignment A2.2.
  \item[$(C3)$] Dominator gets at least $3a-2b+2x_1+(3d-2)x_2$ points due to the
  assignment A2.3.
  \item[$(C4)$] Dominator gets at least $2a+(2d-2)x_3$ points due to the
  assignment A2.4.
  \item[$(C5)$] Dominator gets at least $1$ point due to the
  assignment A2.4.
  \etmz

 Thus, the game starts on $G^*=G_1$ with
 $p(G_1)=a\cdot n$, and ends with a residual graph whose value equals 0.
 Recall that Dominator plays a purely greedy strategy.
  Our goal is to prove that the
 average  decrease in $p(G)$ is at least $s$ points in a  turn.

\paragraph{Phase 1}
  In each turn, the player either selects a W-vertex which turns red
  and hence $p(G)$ decreases by at least $a$ points; or he
  selects a B-vertex $v$ which has a W-neighbor $u$. In the latter
  case the changes $v$: B$\rightarrow$R and $u$: W$\rightarrow$B/R
  together yield a decrease of at least $b+(a-b)=a$ points.
  Therefore, Staller gets at least $a$ points in each of his turns
  in Phase 1. By condition $(C1)$, Dominator seizes at least
  $5a-4b$ points and therefore, in any two consecutive turns $p(G)$
  decreases by at least $6a-4b=2s$ points. As Dominator starts,  the
  following statement follows.
  \bl \label{lemma-4-ph1}
  In   Phase 1, the average decrease of $p(G)$ in a turn is at
  least $s$ points.
  \el
  Concerning the structure of the residual graph obtained at the end
  of this phase, we prove the following properties.
  \bl \label{lemma-4-str1} At the end of Phase 1,
  \tmz
  \item[$(i)$] If\/ $v$ is a W-vertex, then\/ $d^W(v) \le 3$.
  \item[$(ii)$] If\/ $v$ is a B-vertex, then\/ $d(v) \le 4$.
  \etmz \el
  \pf At the end of the phase, we have a residual graph $G_i$ in which Dominator cannot
  get $5a-4b$ or more points.
  Assuming a W-vertex $v$ with  W-neighbors $u_1$, $u_2$, $u_3$ and
  $u_4$, Dominator could play $v$ and the changes
  $v$: W$\rightarrow$R and $u_1, u_2, u_3, u_4$: W$\rightarrow$B/R
  would result in a decrease of at least $a+4(a-b)=5a-4b$ points,
  which is a contradiction.
  In the other case, the choice of a B-vertex  which has five W-neighbors
  would yield a gain of at least $b+5(a-b)=5a-4b$ points, which is a
  contradiction again.
  \qed

  \bsk

  \paragraph{Phase 2} In this phase we apply the value assignment
  A2.2. By Lemma~\ref{lemma-4-str1}$(ii)$, each B-vertex  has degree
  smaller than or equal to 4. Moreover by the definition of A2.2 and
  by Lemma~\ref{lemma-count}, in the $j$th turn the value of a
  B-vertex $u$ decreases by at least $(d_{G_j}(u)-d_{G_{j+1}}(u))
  x_1 $ points.
  \bl \label{lemma-4-ph2}
  In  Phase 2, the average decrease of $p(G)$  in a turn is at
  least $s$ points.
  \el
   \pf
   If a W-vertex $v$ is played, each of its neighbors has a
   decrease of at least $x_1$ points in its value, no matter whether this
   change on the neighbor is
    B$_i\rightarrow$B$_{i-1}$ or   B$_1\rightarrow$R or
    W$\rightarrow$B$_i$/R. Then, playing a W-vertex results in at least
    $a+d\cdot x_1$ point decrease in $p(G)$.

    In the other case, when the played vertex $v$ is blue, the
    decrease in its value  is at least $b-3x_1$. As $v$ has a
    W-neighbor $u$, whose W-degree is at most 3, the change
    $u$: W$\rightarrow$B$_i$/R ($i\le 3$) yields further at least
        $a-(b-x_1)$ points gain; and since $u$ has at least $d-4$ B-neighbors different
    from $v$, the total decrease in $p(G)$ is at least
    $(b-3x_1)+a-(b-x_1)+(d-4)x_1=a+(d-6)x_1$. This yields that
    Staller gets at least $a+(d-6)x_1$ points whenever a white or a blue
    vertex is played by him.

    Complying with $(C2)$, each move of Dominator results in a gain
    of at least $4a-3b+(4d-6)x_1$ and consequently, in any two
    consecutive turns of Phase 2,  $p(G)$ decreases by at least
    $5a-3b+(5d-12)x_1=2s$ points. This proves the lemma. \qed

    \bl \label{lemma-4-str2} At the end of Phase 2,
  \tmz
  \item[$(i)$] If\/ $v$ is a W-vertex, then\/ $d^W(v) \le 2$.
  \item[$(ii)$] If\/ $v$ is a B-vertex, then\/ $d(v) \le 3$.
  \etmz \el
  \pf To prove $(i)$, assume that  Dominator selects  a W-vertex $v$ with W-neighbors $u_1$,
  $u_2$ and $u_3$. Remark that each $u_\ell$ may have at most two
  W-neighbors different from $v$.
  Therefore, the changes
  $v$: W$\rightarrow$R and $u_1, u_2, u_3$: W$\rightarrow$B$_i$/R ($i\le 2$)
  give at least $a+3(a-b+2x_1)$ points to Dominator.
  In addition, each of $v, u_1,u_2$ and $u_3$  has at least $d-3$
  B-neighbors.  Hence the total decrease in $p(G)$ is at least
  $a+3(a-b+2x_1)+4(d-3)x_1=4a-3b+(4d-6)x_1$ points.
  In this case, Dominator's turn would belong to Phase 2. Hence for
  every W-vertex $v$, $d^W(v) \le 2$ must hold at the end of Phase
  2.

  Part $(ii)$ can be shown in a similar way but here we can refer to the property $(i)$ proved above.
  The selection of a
  B-vertex $v$ which has four  W-neighbors, say  $u_1, u_2, u_3$ and
  $u_4$, would cause the color changes
  $v$: B$_4\rightarrow$R and $u_1, u_2, u_3, u_4$: W$\rightarrow$B$_i$/R (where $i\le 2$, due to part
  $(i)$). Moreover each $u_j$ has at least $d-3$ B-neighbors
  different from $v$. These would give a gain of at least
  $$b+4(a-b+2x_1)+4(d-3)x_1=4a-3b+(4d-4)x_1>4a-3b+(4d-6)x_1$$
  point   to Dominator, which is impossible at the end of Phase
  2. This verifies part $(ii)$.
  \qed

  \bsk

  \paragraph{Phase 3} Here we apply the value assignment
  A2.3. By Lemma~\ref{lemma-4-str2}$(ii)$, each B-vertex $v$ has degree
  $d(v) \le 3$ and hence, A2.3 is defined for all vertices of the
  residual graph. We observe concerning this phase that whenever the
  degree of a   B-vertex $v$ is reduced by $y$, its value  decreases
  by at least $y x_2$ points.
 \bl \label{lemma-4-ph3}
 In   Phase 3, the average decrease of $p(G)$ in a turn is at
  least $s$ points. \el
  \pf If Staller plays a W-vertex, he gets at least $a+dx_2$ points.
  In the other case,  he plays a B-vertex $v$ which has a
  W-neighbor $u$. By Lemma~\ref{lemma-4-str2}, $d^W(u) \le 2$ and
  hence,  $u$ has at least $d-3$ B-neighbors different from $v$.
  The changes
  $v$: B$_i\rightarrow$R and $u$: W$\rightarrow$B$_i$/R (where $i\le
  2$), together with the changes on the further B-neighbors of $u$,
  yields a decrease of at least
  $$(b-x_1-2x_2)+ a-(b-x_1-x_2)+(d-3)x_2=a+(d-4)x_2$$
  in $p(G)$. Therefore, Staller gets at least $a+(d-4)x_2$ points in
  each of his turns.
  By condition $(C2)$, we have a lower bound on the gain of
  Dominator as well.  These yield the sum
  $$4a-2b+2x_1+(4d-6)x_2=2s$$
  for any two consecutive turns of Phase 3, and we can conclude that
  the average is at least $s$ indeed.  \qed
  \bsk

  \bl \label{lemma-4-str3} At the end of Phase 3,
  \tmz
  \item[$(i)$] If\/ $v$ is a W-vertex, then\/ $d^W(v) \le 1$.
  \item[$(ii)$] If\/ $v$ is a B-vertex, then\/ $d(v) \le 2$.
  \etmz \el
  \pf At the end of Phase 3 we have a residual graph $G_i$,
  in which the choice of any vertex decreases $p(G_i)$ by strictly
  less than $3a-2b+2x_1+(3d-2)x_2$ points.

  \nin $(i)$ Playing a W-vertex $v$, which has two W-neighbors
  say $u_1$ and $u_2$, results in the changes
  $v$: W$\rightarrow$R and $u_1, u_2$: W$\rightarrow$B$_1$/R; additionally
  $d^B(v)+d^B(u_1)+d^B(u_2) \ge 3(d-2)$. This means a decrease  of at
  least $$a+2(a-b+x_1+2x_2)+3(d-2)x_2=3a-2b++2x_1+(3d-2)x_2$$
  in $p(G_i)$. This cannot be the case; so  each W-vertex has either zero or exactly one W-neighbor in
  $G_i$.

  \nin $(ii)$ Now suppose that a B-vertex $v$ with W-neighbors $u_1$, $u_2$ and $u_3$ is
  played in $G_i$. We have already seen that $d^W(u_j) \le 1$ holds
  for every W-vertex $u_j$ in $G_i$. Then, we have the changes
  $v$: B$_3\rightarrow$R and $u_1, u_2, u_3$: W$\rightarrow$B$_1$/R.
  Further, each vertex from
  $\{u_1, u_2, u_3\}$ has at least $d-2$ B-neighbors different from
  $v$. Hence, the total gain of the player would be at least
  $$b-x_1+3(a-b+x_1+2x_2)+3(d-2)x_2>3a-2b++2x_1+(3d-2)x_2.$$
  This contradiction proves $(ii)$. \qed
  \bsk

  \bsk

  \paragraph{Phase 4} First, we change to assignment
  A2.3. By Lemma~\ref{lemma-4-str3}$(ii)$, in any residual graph of
  Phase 4,
  the W-vertices induce a subgraph consisting
   of isolated vertices and $P_2$-components; moreover, each blue
   vertex has at most 2 (white) neighbors.
   Moreover, by Table~\ref{table-deg4} and Lemma~\ref{counting}, if
   a B-vertex loses $y$ W-neighbors in a turn, its value is reduced
   by at least $yx_3$ points.

 \bl \label{lemma-4-ph4}
 In   Phase 4, the average decrease of $p(G)$ in a turn is at
  least $s$ points. \el
  \pf If Staller selects a W-vertex $v$,  each neighbor $u$ of $v$ has a
  decrease of at least $x_3$ in its value. Hence, the total decrease
  in $p(G)$ is not smaller than
  $b-x_1-x_2+dx_3$.

  If Staller selects a B-vertex $v$, the change is either
  $v$: B$_2\rightarrow$R or $v$: B$_1\rightarrow$R, it means at
  least $(b-x_1-x_2-x_3)$-point gain. As $d(v)\ge 1$, we necessarily
  have a W-neighbor $u$ of $v$ whose change is
  $u$: W$\rightarrow$B$_1$/R. Further, $u$ has at least $d-2$
  B-neighbors different from $v$.
  Therefore, the decrease in $p(G)$ is at least
  $$(b-x_1-x_2-x_3)+a-(b-x_1-x_2-x_3)+(d-2)x_3=a+(d-2)x_3.
  $$
  Hence, in any case, Staller gets at least $a+(d-2)x_3$ points in a
  turn of his own. By $(C4)$, Dominator gets at least $2a+(2d-2)x_3$
  point in each of his turns and as follows, the average gain is at
  least
  $$\frac{1}{2}(a+(d-2)x_3+2a+(2d-2)x_3=s$$
  points as stated.  \qed
  \bsk

  \bl \label{lemma-4-str4} At the end of Phase 4,
  \tmz
  \item[$(i)$] Every  W-vertex has only B-neighbors.
  \item[$(ii)$] Every B-vertex has exactly one W-neighbor.
  \etmz \el
  \pf Consider $G_i$ which is the residual graph obtained at the end
  of Phase 4. As $(C4)$ is not true, Dominator cannot get
  $2a+(2d-2)x_3$ or more points in the $i$th turn.
  By Lemma~\ref{lemma-4-str3}, if $(i)$ is not true, we have a
  "white-pair" $(v,u)$, where $u$ is the only W-neighbor of $v$ and
  vice versa. Then, the choice of $v$ would result the changes
  $v,u$: W$\rightarrow$R. This, together with the fact $d^B(v)+d^B(u) \ge
  2d-2$, implies that the decrease in $p(G_i)$ is at least
  $2a+(2d-2)x_3$, which is a contradiction. Thus, $(i)$ is true.

  To prove $(ii)$ we suppose for a contradiction that a B-vertex $v$
  has two W-neighbors $u_1$ and $u_2$. By $(i)$, these neighbors are
  "single-white" vertices and they
  turn to red if $v$ is played; in addition both $u_1$ and $u_2$ has
  at least $d-1$ B-neighbors different from $v$.
  Hence, selecting $v$ Dominator could seize at least
  $$(b-x_1-x_2)+2a+2(d-1)x_3 > 2a+(2d-2)x_3$$
  points, which is a contradiction again.
  \qed
  \bsk

  \paragraph{Phase 5}
   By Lemma~\ref{lemma-4-str4}$(ii)$, the residual graphs occurring
   in this phase have simple structure, each of their components is
   a star of order at least $d+1$ whose center is white and the
   leaves are blue. Then, in each turn of Phase 5 exactly one such a star component becomes completely red,
   no matter whether a white or a blue vertex is played.
   Then, the value of the residual graph is decreased by at least
   $a+d(b-x_1-x_2-x_3)=s$ points in each single turn.

 \bl \label{lemma-4-ph5}
 In   Phase 5, the average decrease of $p(G)$ in a turn is at
  least $s$ points. \el

  By Lemmas \ref{lemma-4-ph1}, \ref{lemma-4-ph2}, \ref{lemma-4-ph3},
  \ref{lemma-4-ph4} and \ref{lemma-4-ph5}, the average decrease per
  turn in the residual graph is at least $s$ for the entire game. As
  $p(G_1)=an$ and the changes between assignments  nowhere caused
  increase in $p(G)$, the domination game where Dominator plays
  the described greedy strategy yields a game with at most $an/s$
  turns. This establishes Theorem~\ref{main}. \qed

  \bsk

\section{Concluding remarks on the Staller-start game}

In our main theorem, we do not give upper bound on $\gamma_g'(G)$
for graphs with $\delta(G)\ge d\ge 4$. It is quite clear from the
proof that we can establish the same upper bound on $\gamma_g'(G)$
as proved for $\gamma_g(G)$. Moreover, a slight improvement on it is
also possible. We close the paper with this complicated formula.
\bsk

If Staller begins the game, we index this starting turn by zero and
take it into Phase 0. Then, from the first turn of Dominator, it
continues as in the proof of Theorem~\ref{main}. In the turn of
Phase 0, Staller gets at least $$a+d(a-b)=s+30d^5  -227d^4+ 637d^3
-786d^2+ 360d$$ points. In later phases, the average decrease
remains at least $s$. This proves that

$$\gamma_g'(G)  \le
        \frac{(30d^4-56d^3-258d^2+708d-432)n-30d^5  +227d^4    -637d^3 +786d^2-
        360d}{90d^4-390d^3+348d^2+348d-432}$$
        holds for every $d\ge 4$ and for every graph $G$ of minimum degree not smaller then
        $d$.


\end{document}